\documentclass[12pt,a4paper]{article}
\usepackage{amssymb}
\usepackage{amsmath}
\usepackage{amsthm}
\usepackage[arrow,matrix]{xy}
\let\mod=\undefined
\DeclareMathOperator{\coker}{Coker}

\DeclareMathOperator{\Ext}{Ext}
\DeclareMathOperator{\GL}{GL}
\DeclareMathOperator{\Hom}{Hom}
\DeclareMathOperator{\im}{Im}

\DeclareMathOperator{\Ker}{Ker}
\DeclareMathOperator{\mod}{mod}

\DeclareMathOperator{\rk}{rk}
\DeclareMathOperator{\Reg}{Reg}

\DeclareMathOperator{\Sing}{Sing}

\newcommand{\BB}{{\mathbb B}}
\newcommand{\BM}{{\mathbb M}}

\newcommand{\BP}{{\mathbb P}}
\newcommand{\BZ}{{\mathbb Z}}
\newcommand{\CA}{{\mathcal A}}
\newcommand{\CB}{{\mathcal B}}
\newcommand{\CC}{{\mathcal C}}
\newcommand{\CE}{{\mathcal E}}
\newcommand{\CF}{{\mathcal F}}
\newcommand{\CG}{{\mathcal G}}
\newcommand{\CH}{{\mathcal H}}
\newcommand{\CO}{{\mathcal O}}
\newcommand{\CQ}{{\mathcal Q}}
\newcommand{\CS}{{\mathcal S}}

\newcommand{\CU}{{\mathcal U}}
\newcommand{\CV}{{\mathcal V}}
\newcommand{\CW}{{\mathcal W}}
\newcommand{\CX}{{\mathcal X}}
\newcommand{\CY}{{\mathcal Y}}
\newcommand{\CZ}{{\mathcal Z}}

\newcommand{\ov}{\overline}

\newcommand{\bsmatrix}[1]{\left[\begin{smallmatrix} #1%
 \end{smallmatrix}\right]}
\newcommand{\psmatrix}[1]{\left(\begin{smallmatrix} #1%
 \end{smallmatrix}\right)}
\newtheorem{thm}{Theorem}[section]
\newtheorem{cor}[thm]{Corollary}
\newtheorem{lem}[thm]{Lemma}
\newtheorem{prop}[thm]{Proposition}

\numberwithin{equation}{section}
\begin{document}
\title{Singularities of orbit closures in module varieties and cones 
 over rational normal curves
 \footnotetext{Mathematics Subject Classification (2000): %
14L30 (Primary); 14B05, 16G10, 16G20 (Secondary).}
 \footnotetext{Key Words and Phrases: Module varieties,
 orbit closures, types of singularities, rational normal curves.}}
\author{Grzegorz Zwara}
%
\maketitle

\begin{abstract}
Let $N$ be a point of an orbit closure $\ov{\CO}_M$ in a module variety
such that its orbit $\CO_N$ has codimension two in $\ov{\CO}_M$.
We show that under some additional conditions the pointed variety 
$(\ov{\CO}_M,N)$ is smoothly equivalent to a cone over a rational normal
curve.  
\end{abstract}

\section{Introduction and the main results}

Throughout the paper, $k$ denotes an algebraically closed field,
$A$ denotes a finitely generated associative $k$-algebra with
identity, and by a module we mean a left $A$-module whose
underlying vector space is finite dimensional.
Let $d$ be a positive integer and denote by $\BM(d)$
the algebra of $d\times d$-matrices with coefficients in $k$.
The set $\mod_A(d)$ of algebra homomorphisms
$A\to\BM(d)$ has a natural structure of an affine variety.
Indeed, if $A\simeq k\langle X_1,\ldots,X_t\rangle/I$ for some
two-sided ideal $I$, then $\mod_A(d)$ can be identified with
the closed subset of $(\BM(d))^t$ given by vanishing of the
entries of all matrices $\rho(X_1,\ldots,X_t)$, $\rho\in I$.
Moreover, the general linear group $\GL(d)$ acts on $\mod_A(d)$
by conjugations
$$
g\star(M_1,\ldots,M_t)=(gM_1g^{-1},\ldots,gM_tg^{-1}),
$$
and the $\GL(d)$-orbits in $\mod_A(d)$ correspond bijectively
to the isomorphism classes of $d$-dimensional modules.
We shall denote by $\CO_M$ the $\GL(d)$-orbit in $\mod_A(d)$
corresponding to a $d$-dimensional module $M$.
An interesting problem is to study geometric properties of the
Zariski closure $\ov{\CO}_M$ of an orbit $\CO_M$ in $\mod_A(d)$.
We refer to \cite{BB}, \cite{BZ2}, \cite{Bmin}, \cite{Bext}, 
\cite{Zuni}, \cite{Zcodim1} and \cite{Zcodim2} for some results 
in this direction.

Following Hesselink (see~\cite[(1.7)]{Hes})
we call two pointed varieties $(\CX,x_0)$ and $(\CY,y_0)$
smoothly equivalent if there are smooth morphisms
$f:\CZ\to\CX$, $g:\CZ\to\CY$ and a point $z_0\in\CZ$
with $f(z_0)=x_0$ and $g(z_0)=y_0$.
This is an equivalence relation and the equivalence classes
will be denoted by $\Sing(\CX,x_0)$ and called the types
of singularities.
Obviously the regular points of the varieties form one type
of singularity, which we denote by $\Reg$.
Let $M$ and $N$ be $d$-dimensional modules such that $N$ is 
a degeneration of $M$, i.e., $\CO_N\subseteq\ov{\CO}_M$.
We shall write $\Sing(M,N)$ for $\Sing(\ov{\CO}_M,n)$, where $n$
is an arbitrary point of $\CO_N$.
By \cite[Thm.1.1]{Zcodim1}, $\Sing(M,N)=\Reg$ provided 
$\dim\CO_M-\dim\CO_N=1$.

Assume now that $\dim\CO_M-\dim\CO_N=2$.
It was shown recently (\cite[Thm.1.3]{Zcodim2}) that $\Sing(M,N)=\Reg$ 
provided $A$ is the path algebra of a Dynkin quiver.
However, other types of singularities occur if $A$ is the path 
algebra of the Kronecker quiver (see \cite{BB}); for instance, 
the type of an isolated singularity of the surface
$$
\CC_r=\left\{(x_0,\ldots,x_r)\in k^{r+1};\;x_ix_j=x_lx_m\text{ if }
i+j=l+m\right\},\quad r\geq 1.
$$ 
The surface $\CC_r$ is the affine cone over the rational normal curve
(or Veronese curve) of degree $r$ (see \cite{Harr} for an introduction 
to rational normal curves).
Obviously $\Sing(\CC_r,0)=\Reg$ if and only if $r=1$.
The main theorem of the paper shows that the above types of singularities
occur quite frequently for orbit closures in module varieties.
From now on, we abbreviate $\dim_k\Hom_A(U,V)$ by $[U,V]$ for any modules
$U$ and $V$.

\begin{thm} \label{mainthm}
Let $M$ and $N$ be modules such that $\dim\CO_M-\dim\CO_N=2$,
$\CO_N\subseteq\ov{\CO}_M$ and $N=U\oplus V$ for some modules $U$ and $V$
satisfying
$$
[U,M]=[U,N]\qquad\text{and}\qquad[M,V]=[N,V].
$$
Then $\Sing(M,N)=\Sing(\CC_r,0)$ for some $r\geq 1$.
\end{thm}

As an application of this theorem we get the following result about
preprojective modules (see Section~\ref{proof2} for the definition).

\begin{thm} \label{mainthm2}
Assume that the algebra $A$ is finite dimensional.
Let $M$ and $N$ be preprojective modules such that 
$\CO_N\subseteq\ov{\CO}_M$ and $\dim\CO_M-\dim\CO_N=2$.
Then $\Sing(M,N)=\Sing(\CC_r,0)$ for some $r\geq 1$.
\end{thm}

One says that a module $N$ is a minimal degeneration of a module $M$ 
provided $\ov{\CO}_N\varsubsetneq\ov{\CO}_M$ but there is no module 
$L$ with $\ov{\CO}_N\varsubsetneq\ov{\CO}_L\varsubsetneq\ov{\CO}_M$.
Moreover, two modules are called disjoint if they have no nonzero 
direct summand in common. 
Combining the main results in \cite{Bgeo}, \cite{BF} and \cite{Zcodim1} 
with Theorem~\ref{mainthm2}, we get the following consequence.

\begin{cor}
Let $A$ be the path algebra of an extended Dynkin quiver.
Let $M$ and $N$ be two disjoint preprojective modules such that 
$N$ is a minimal degeneration of $M$.
Then $\Sing(M,N)=\Sing(\CC_r,0)$ for some integer $r\geq 1$.
\end{cor}

We consider in Section~\ref{degen} some properties of short exact 
sequences, dimensions of homomorphism spaces and degenerations 
of modules.
Sections~\ref{genquot}, \ref{monext} and \ref{extnormal} contain
preliminary results necessary to prove Proposition~\ref{keyprop}, 
which plays a central role in the paper.
Sections~\ref{proof1} and \ref{proof2} are devoted to the proofs 
of Theorems~\ref{mainthm} and \ref{mainthm2}, respectively.
In Section~\ref{example} we give an application of Theorem~\ref{mainthm} 
different from the one described in Theorem~\ref{mainthm2}.  

For a basic background on the representation theory of algebras and
quivers we refer to \cite{ARS} and \cite{Rin}.
The author gratefully acknowledges support from the Polish
Scientific Grant KBN No.\ 1 P03A 018 27.

\section{Degenerations and extensions of modules}
\label{degen}

Let $W$ be a vector space. 
We denote by $\BP W$ the projective space of one-dimensional subspaces in $W$.
Thus we have the quotient map
$$
\pi:\;W\setminus\{0\}\to\BP W,\qquad w\mapsto k\cdot w,
$$ 
of $W\setminus\{0\}$ by $k^\ast$ acting by scalar multiplication.
Moreover, we denote by $\BP S$ the image $\pi(S\setminus\{0\})$ 
for any $k^\ast$-invariant subset $S$ of $W$.

Throughout the paper, we consider short exact sequences of modules 
\begin{equation} \label{sigma}
\sigma:\;0\to U\xrightarrow{f}M\xrightarrow{g}V\to 0,
\end{equation}
where $U$, $M$ and $V$ are modules, $f$ is an $A$-module monomorphism, 
$g$ is an $A$-module epimorphism and $\im(f)=\Ker(g)$.
Any such sequence $\sigma$ determines an element of $\Ext^1_A(V,U)$ 
denoted by $[\sigma]$, such that $[\sigma]=0$ if and only if the sequence 
$\sigma$ splits, and  $[\sigma]=[\tau]$ for another short exact sequence 
$$
\tau:\;0\to U\xrightarrow{f'}W\xrightarrow{g'}V\to 0
$$ 
if and only if there is an $A$-isomorphism $h:M\to W$ making the following
diagram commutative:
$$
\xymatrix{
0\ar[r]&U\ar[r]^f\ar@{=}[d]&M\ar[r]^g\ar[d]^h&V\ar[r]\ar@{=}[d]&0\\
0\ar[r]&U\ar[r]^{f'}&W\ar[r]^{g'}&V\ar[r]&0.
}
$$
We denote by $\CE_M(V,U)$ the subset of $\Ext^1_A(V,U)$ consisting 
of all elements $[\sigma]$, such that the middle term of the exact 
sequence $\sigma$ is isomorphic to $M$.

\begin{lem}
$\CE_M(V,U)$ is a locally closed $k^\ast$-invariant subset
of the vector space $\Ext^1_A(V,U)$.
\end{lem}

\begin{proof}
We use here an interpretation of the vector space $\Ext^1_A(V,U)$ as 
the quotient of the vector space $\BZ^1_A(V,U)$ of cocycles, i.e., the
$k$-linear maps $Z:A\to\Hom_k(V,U)$ satisfying
$$
Z(aa')=Z(a)V(a')+U(a)Z(a'),\qquad\text{for all }a,a'\in A,
$$
by the subspace of coboundaries
$$
\BB^1_A(V,U)=\{hV-Uh;\;h\in\Hom_k(V,U)\}.
$$
Let $d=\dim_k(U\oplus V)$.
We take a subspace $\CW$ of $\BZ^1_A(V,U)$ complementary to $\BB^1_A(V,U)$.
Then the subset $\CW'$ of $\mod_A(d)$ consisting of the algebra homomorphism
of the block form $\begin{bmatrix}U&Z\\ 0&V\end{bmatrix}$, where $Z\in\CW$,
is closed.
Hence the claim follows from the fact that $\CE_M(V,U)$ is isomorphic to
the intersection $\CW'\cap\CO_M$. 
\end{proof}

The short exact sequence \eqref{sigma} induces two long exact sequences
of vector spaces
\begin{multline*}
0\to\Hom_A(V,X)\xrightarrow{\Hom_A(g,X)}\Hom_A(M,X)\xrightarrow{\Hom_A(f,X)}
 \Hom_A(U,X)\xrightarrow{\partial}\\
\mbox{}\xrightarrow{\partial}\Ext^1_A(V,X)\xrightarrow{\Ext^1_A(g,X)}
 \Ext^1_A(M,X)\xrightarrow{\Ext^1_A(f,X)}\Ext^1_A(U,X),
\end{multline*}
\begin{multline*}
0\to\Hom_A(X,U)\xrightarrow{\Hom_A(X,f)}\Hom_A(X,M)\xrightarrow{\Hom_A(X,g)}
 \Hom_A(X,V)\xrightarrow{\partial'}\\
\mbox{}\xrightarrow{\partial'}\Ext^1_A(X,U)\xrightarrow{\Ext^1_A(X,f)}
 \Ext^1_A(X,M)\xrightarrow{\Ext^1_A(X,g)}\Ext^1_A(X,V),
\end{multline*}
for any module $X$.
Here, 
$$
\partial(h)=\Ext^1_A(V,h)([\sigma])\qquad\text{and}\qquad 
\partial'(h')=\Ext^1_A(h',U)([\sigma])
$$ 
for any homomorphisms $h:U\to X$ and $h':X\to V$.
These induced exact sequences suggest us to define the following 
two non-negative integers
\begin{align*}
\delta_\sigma(X)&=[U\oplus V,X]-[M,X]=\dim_k\im(\partial)\geq 0,\\
\delta'_\sigma(X)&=[X,U\oplus V]-[X,M]=\dim_k\im(\partial')\geq 0.
\end{align*}
Observe that the sequence $\sigma$ splits if and only if
$$
M\simeq U\oplus V \iff \delta_\sigma(U)=0 \iff \delta'_\sigma(V)=0
$$
(see \cite[Lem.2.3]{Zcodim2}).
Taking the pushout of $\sigma$ under $h$ leads to the following 
commutative diagram with exact rows:
$$
\xymatrix{
\sigma:&0\ar[r]&U\ar[r]^f\ar[d]_h&M\ar[r]^g\ar[d]_j&V\ar[r]\ar@{=}[d]&0\\
\sigma':&0\ar[r]&X\ar[r]^{f'}&W\ar[r]^{g'}&V\ar[r]&0.
}
$$
Then $[\sigma']=\Ext^1_A(V,h)([\sigma])=\partial(h)$.
Moreover, the short sequence
$$
\tau:\quad 0\to U\xrightarrow{\psmatrix{f\\ h}}M\oplus X
 \xrightarrow{(j,-f')}W\to 0
$$
is exact and
$$
\delta_\sigma(Y)=\delta_{\sigma'}(Y)+\delta_\tau(Y),\qquad
\delta'_\sigma(Y)=\delta'_{\sigma'}(Y)+\delta'_\tau(Y),
$$
for any module $Y$.
In particular, we obtain the following fact.

\begin{cor} \label{pushout}
Let $\sigma'$ be a pushout of a short exact sequence $\sigma$.
Then $\delta_{\sigma'}(Y)\leq\delta_\sigma(Y)$ and 
$\delta'_{\sigma'}(Y)\leq\delta'_\sigma(Y)$ for any module $Y$.
\end{cor}

We shall need the following technical result.
 
\begin{lem} \label{cancelinseq}
Let $\sigma:\;0\to U\oplus Y^i\xrightarrow{f}W\to V\to 0$ be a short 
exact sequence such that $\delta_\sigma(Y)<i$. 
Then $W$ is isomorphic to $W'\oplus Y$ and there is a short exact 
sequence
$$
\sigma':\;0\to U\oplus Y^{i-1}\to W'\to V\to 0
$$
for some module $W'$.
Moreover, $\delta_{\sigma'}(X)=\delta_\sigma(X)$ for any module $X$.
\end{lem}

\begin{proof}
We consider the exact sequence
$$
0\to\Hom_A(V,Y)\to\Hom_A(W,Y)\xrightarrow{\Hom_A(f,Y)}\Hom_A(U\oplus Y^i,Y)
$$
induced by $\sigma$.
Let $\CH$ be the $i$-dimensional subspace of $\Hom_A(U\oplus Y^i,Y)$ spanned 
by $i$ canonical projections.
Since the codimension of the image of $\Hom_A(f,Y)$ in $\Hom_A(U\oplus Y^i,Y)$
equals $\delta_\sigma(Y)<\dim_k\CH$, the image contains a nonzero element
$h\in\CH$.
Thus $h$ factors through $f$.
We decompose 
$$
f=(f_0,f_1,\ldots,f_i)\qquad\text{and}\qquad 
h=(0,\lambda_1\cdot 1_Y,\ldots,\lambda_i\cdot 1_Y)\neq 0.
$$
Then $\lambda_j\neq 0$ and consequently $f_j$ is a section, for some 
$j\in\{1,\ldots,i\}$.
The remaining part of the proof is now straightforward (for instance, use
a dual version of \cite[Lem.2.4]{Zcodim2}).
\end{proof}

The next result (\cite{Zgiv}) gives a characterization of degenerations 
of modules.

\begin{thm} \label{specialseq}
Let $M$ and $N$ be modules.
Then the inclusion $\CO_N\subseteq\ov{\CO}_M$ is equivalent
to each of the following conditions:
\begin{enumerate}
\item[(1)] There is a short exact sequence 
 $0\to Z\to Z\oplus M\to N\to 0$ for some module $Z$.
\item[(2)] There is a short exact sequence
 $0\to N\to M\oplus Z'\to Z'\to 0$ for some module $Z'$.
\end{enumerate}
\end{thm}

We shall need the following two well known direct consequences.

\begin{cor} \label{shortdeg}
Let $\sigma:\;0\to U\to M\to V\to 0$ be a short exact sequence.
Then $\CO_{U\oplus V}\subseteq\ov{\CO}_M$.
\end{cor}

\begin{proof}
We apply Theorem~\ref{specialseq} to a direct sum of $\sigma$ and
the short exact sequence $0\to 0\to U\xrightarrow{1_U}U\to 0$.
\end{proof}

\begin{cor} \label{deghom}
Let $M$ and $N$ be modules such that $\CO_N\subseteq\ov{\CO}_M$.
Then 
$$
[M,Y]\leq[N,Y]\qquad\text{and}\qquad[Y,M]\leq[Y,N]
$$
for any module $Y$.
\end{cor}

\begin{proof}
Observe that 
$$
\delta_\sigma(Y)=[N,Y]-[M,Y]\qquad\text{and}\qquad
\delta'_\sigma(Y)=[Y,N]-[Y,M]
$$
for any short exact sequence $\sigma$ of the form
$0\to Z\to Z\oplus M\to N\to 0$.
\end{proof}

\section{Generic quotients}
\label{genquot}

Throughout the section, $U$, $M$ and $V$ are modules such that 
there is a short exact sequence
$$
\sigma:\;0\to U\to M\to V\to 0
$$
with $\delta_\sigma(V)=0$.
Let $d=\dim_kV$ and $\CQ(U,M)$ be the subset of $\mod_A(d)$ consisting 
of all points corresponding to quotients of $M$ by $U$.
Obviously, the set $\CQ(U,M)$ is $\GL(d)$-invariant.
Moreover, it is constructible and irreducible, by \cite[2.3]{Bext}.

\begin{lem} \label{Vgeneric}
The orbit $\CO_V$ is an open subset of $\ov{\CQ(U,M)}$.
\end{lem}

\begin{proof}
Let $\CA$ be the closed subset of $\ov{\CQ(U,M)}\times\Hom_k(k^d,V)$ 
consisting of all pairs $(V',h)$ such that $h$ belongs to $\Hom_A(V',V)$.
Then the projection $\pi:\CA\to\ov{\CQ(U,M)}$ is a surjective regular 
morphism.
We know that $U\oplus V'$ belongs to $\ov{\CO}_M$
for any $V'\in\CQ(U,M)$, by Corollary~\ref{shortdeg}.
Hence the same holds for arbitrary $V'\in\ov{\CQ(U,M)}$.
Consequently, 
$$
[U\oplus V',V]\geq[M,V]=[U\oplus V,V]\quad\text{and}\quad
[V',V]\geq[V,V],
$$
by Corollary~\ref{deghom} and since $\delta_\sigma(V)=0$.
We consider the closed subset $\CA'$ of $\CA$ consisting of the pairs 
$(V',h)$ such that $h$ is not invertible.
Observe that 
$$
\dim(\pi|_{\CA'})^{-1}(V')=
\begin{cases}[V,V]-1,&\text{if }V'\simeq V,\\ 
 [V',V],&\text{if }V'\not\simeq V.\end{cases}
$$
Hence $(\pi|_{\CA'})^{-1}(\CO_V)$ is an open subset of $\CA'$.
Applying the inverse image of the section
$$
\ov{\CQ(U,M)}\to\CA',\quad V'\mapsto(V',0),
$$
we get that $\CO_V$ is an open subset of $\ov{\CQ(U,M)}$. 
\end{proof}

Let $\CS_V(U,M)$ denote the set of all $A$-monomorphisms $f:U\to M$ 
with $\coker(f)\simeq V$.

\begin{lem} \label{Sopen}
$\CS_V(U,M)$ is an open subset of $\Hom_A(U,M)$.
\end{lem}

\begin{proof}
We use a construction described in \cite[2.1]{Bmin}.
Let $\CG$ be the subset of $\Hom_A(U,M)\times\Hom_k(k^d,M)$ consisting of 
pairs $g=(g_1,g_2)$ such that the map $g:U\oplus k^d\to M$ is invertible.
Then the projection $\pi:\CG\to\Hom_A(U,M)$ is a composition 
of an open immersion followed by a trivial vector bundle.
Moreover, if $g\in\CG$ then $g^{-1}\star M$ has a block form
$$
\begin{bmatrix}U&Z\\ 0&V'\end{bmatrix},
$$
where $V'$ belongs to the set $\CQ(U,M)$.
This leads to a regular morphism of varieties 
$$
\theta:\,\CG\to\ov{\CQ(U,M)},\quad g\mapsto V'.
$$
Then $\theta^{-1}(\CO_V)$ is open in $\CG$, by Lemma~\ref{Vgeneric}.
Since the map $\pi$ is open, the set $\CS_V(U,M)=\pi(\theta^{-1}(\CO_V))$ 
is open in $\Hom_A(U,M)$.
\end{proof}

\section{Monomorphisms and extensions}
\label{monext}

Throughout the section, $U$, $M$ and $V$ are modules such that 
there is a short exact sequence
$$
\sigma:\;0\to U\xrightarrow{f}M\to V\to 0
$$
with $\delta_\sigma(V)=0$ and $\delta'_\sigma(V)=1$.

\begin{lem} \label{fromStoP}
Let $\tau:\;0\to U\xrightarrow{f'}W\to V\to 0$ be a nonsplittable short
exact sequences such that $f$ factors through $f'$ (equivalently, 
the pushout of $\tau$ under $f$ is a splittable sequence).
Then $W\simeq M$ and the equality $k\cdot[\tau]=k\cdot[\sigma]$ holds
in $\BP\Ext^1_A(V,U)$.
\end{lem}

\begin{proof}
Since $f$ factors through $f'$, we get the following commutative diagram
$$
\xymatrix{
\tau:&0\ar[r]&U\ar[r]^{f'}\ar@{=}[d]&W\ar[r]\ar[d]&V\ar[r]\ar[d]^h&0\\
\sigma:&0\ar[r]&U\ar[r]^f&M\ar[r]&V\ar[r]&0.
}
$$
We consider the long exact sequence 
$$
0\to\Hom_A(V,U)\to\Hom_A(V,M)\to\Hom_A(V,V)\xrightarrow{\partial}
\Ext^1_A(V,U)
$$
induced by $\sigma$. 
Then $\partial(1_V)=[\sigma]$ and $\partial(h)=[\tau]$.
Now the claim follows from the equalities
$1=\delta'_\sigma(V)=\dim_k\im(\partial)$ and since 
$[\sigma]\neq 0\neq[\tau]$.
\end{proof}

The above lemma shows the existence of a unique surjective map 
$$
\xi_{U,M,V}:\CS_V(U,M)\to\BP\CE_M(V,U),
$$
such that $\xi_{U,M,V}(f')=k\cdot[\sigma']$ for any short exact 
sequence of the form
$$
\sigma':\;0\to U\xrightarrow{f'}M\to V\to 0.
$$

\begin{lem}
$\xi_{U,M,V}$ is a regular morphism of varieties.
\end{lem}

\begin{proof}
Let $d=\dim_kV$.
We introduce three new varieties and four regular morphisms:
$$
\xymatrix{
&\CA\ar[dl]_{\alpha}\ar[dr]^{\varphi}&&\CC\ar[dl]_{\beta}\ar[dr]^{\psi}\\
\CS_V(U,M)&&\CB&&\BP\CE_M(V,U).
}
$$  
Here $\CA$ is the subset of $\CS_V(U,M)\times\Hom_k(k^d,M)$ consisting 
of all pairs $g=(g_1,g_2)$ such that the map $g:U\oplus k^d\to M$ 
is invertible; $\CB$ is the subset of $\CO_M$ consisting of all elements 
having a block form
$$
\begin{bmatrix}U&Z\\ 0&W\end{bmatrix}
$$
with $W\in\CO_V$; and $\CC$ is the subset of $\CB\times\Hom_k(V,k^d)$
consisting of all pairs $\left(\bsmatrix{U&Z\\ 0&W},h\right)$ such that 
$h$ is an isomorphism in $\Hom_A(V,W)$.
Moreover, $\alpha$ and $\beta$ are obvious projections, 
and $\varphi(g)=g^{-1}\star M$.
If $c=\left(\bsmatrix{U&Z\\ 0&W},h\right)$ belongs to $\CC$, then
$$
\begin{bmatrix}1_U&0\\ 0& h^{-1}\end{bmatrix}\star
\begin{bmatrix}U&Z\\ 0&W\end{bmatrix}
=\begin{bmatrix}U&Zh\\ 0&V\end{bmatrix}\in\CO_M,
$$
where $Zh$ belongs to $\BZ^1_A(V,U)\setminus\BB^1_A(V,U)$.
We define $\psi(c)$ as the image of $Zh$ via the following composition 
of canonical quotients
$$
\BZ^1_A(V,U)\setminus\BB^1_A(V,U)\to\Ext^1_A(V,U)\setminus\{0\}
\to\BP\Ext^1_A(V,U).
$$
One can find that $\alpha$ and $\beta$ are surjective maps being compositions 
of open immersions followed by vector bundles (see for instance \cite[2.1]{Bmin}).
Hence $\alpha$ and $\beta$ have enough local sections, 
i.e., there are open neighbourhoods $\CU$ of $s$ and $\CV$ of $b$ together 
with regular morphisms (sections) $\alpha':\CU\to\CA$, $\beta':\CV\to\CC$ 
such that $\alpha\alpha'=1_\CU$ and $\beta\beta'=1_\CV$, for any 
$s\in\CS_V(U,M)$ and $b\in\CB$.  

We take $s\in\CS_V(U,M)$ and consider a local section $\alpha':\CU\to\CA$ 
of $\alpha$ with $\CU$ containing $s$.
Let $a=\alpha'(s)$ and $b=\varphi(a)$.
We take also a neighbourhood $\CV$ of $b$ and a local section $\beta':\CV\to\CC$ 
of $\beta$.
Then $\CU'=(\varphi\alpha')^{-1}(\CV)$ is an open neighbourhood of $s$ and
it follows from our constructions that 
$\xi_{U,M,V}(s')=\psi\beta'\varphi\alpha'(s')$ for any $s'\in\CU'$.
Therefore $\xi_{U,M,V}|_{\CU'}$ is a regular morphism.
\end{proof}

The set $\CS_V(U,M)$ is irreducible, by Lemma~\ref{Sopen}.
Applying the surjective morphism $\xi_{U,M,V}$ we get the following
consequence.

\begin{cor} \label{EMirreducible}
$\BP\CE_M(V,U)$ is an irreducible variety.
\end{cor}

\section{Extensions and rational normal curves}
\label{extnormal}

Throughout the section, 
$$
\sigma_1:\;0\to U_0\xrightarrow{f_1}U_1\xrightarrow{h_1}V\to 0
$$
is a short exact sequence such that
$$
\delta_{\sigma_1}(U_1)=1,\quad\delta_{\sigma_1}(V)=0\quad\text{and}\quad
\delta'_{\sigma_1}(V)=1.
$$
Our aim is to prove that $\BP\ov{\CE_{U_1}(V,U_0)}$ is a rational normal 
curve contained in $\BP\Ext^1_A(V,U_0)$.

We consider the exact sequence
$$
0\to\Hom_A(V,U_1)\to\Hom_A(U_1,U_1)\xrightarrow{\Hom_A(f_1,U_1)}
\Hom_A(U_0,U_1)
$$
induced by $\sigma_1$.
Since $\delta_{\sigma_1}(U_1)>0$ then there is an $A$-homomorphism 
$g_1:U_0\to U_1$ not belonging to the image of $\Hom_A(f_1,U_1)$.
We construct successive pushouts to obtain the following commutative 
diagram with exact rows
\begin{equation} \label{diag1}
\vcenter{\xymatrix{
\sigma_1:&0\ar[r]&U_0\ar[r]^{f_1}\ar[d]_{g_1}&U_1\ar[r]^{h_1}\ar[d]^{g_2}
 &V\ar[r]\ar@{=}[d]&0\\
\sigma_2:&0\ar[r]&U_1\ar[r]^{f_2}\ar[d]_{g_2}&U_2\ar[r]^{h_2}\ar[d]^{g_3}
 &V\ar[r]\ar@{=}[d]&0\\
\sigma_3:&0\ar[r]&U_2\ar[r]^{f_3}\ar[d]_{g_3}&U_3\ar[r]^{h_3}\ar[d]^{g_4}
 &V\ar[r]\ar@{=}[d]&0\\
&&\vdots&\vdots&\vdots
}}
\end{equation}

Applying Corollary~\ref{pushout} we get the following fact.

\begin{cor} \label{delta0to1}
$\delta_{\sigma_i}(V)=0$ and $\delta'_{\sigma_i}(V)\leq 1$ for any 
$i\geq 1$.
\end{cor}

\begin{lem} \label{fatseq}
There is a short exact sequence
$$
0\to U_0\to U_i\to V^i\to 0
$$
for any $i\geq 1$.
\end{lem}

\begin{proof}
We proceed by induction on $i$.
Assume that there is a short exact sequence
$$
0\to U_0\xrightarrow{f'}U_i\xrightarrow{h'}V^i\to 0
$$
for some $i\geq 1$.
Taking the pushout of $\sigma_{i+1}$ under $h'$ we obtain the following 
commutative diagram with exact rows and columns
$$
\xymatrix{
&&0\ar[d]&0\ar[d]\\
&&U_0\ar@{=}[r]\ar[d]_{f'}&U_0\ar[d]\\
\sigma_{i+1}:&0\ar[r]&U_i\ar[r]^{f_{i+1}}\ar[d]_{h'}
 &U_{i+1}\ar[r]^{h_{i+1}}\ar[d]&V\ar[r]\ar@{=}[d]&0\\
\tau:&0\ar[r]&V^i\ar[r]\ar[d]&W\ar[r]\ar[d]&V\ar[r]&0.\\
&&0&0
}
$$
Thus it suffices to show that $W$ is isomorphic to $V^{i+1}$.
Applying Corollaries~\ref{pushout} and \ref{delta0to1}, we get
$\delta_\tau(V)\leq\delta_{\sigma_{i+1}}(V)=0$.
Hence $\delta_\tau(V^i)=0$.
Therefore the short exact sequence $\tau$ splits and $W$ is isomorphic 
to $V^i\oplus V$.
\end{proof}

\begin{lem}
The short exact sequence $\sigma_i$ splits for some $i$.
\end{lem}

\begin{proof}
Let $i$ be an integer not smaller than $j=\dim_k\Ext^1_A(V,U_0)$.
We get from Lemma~\ref{fatseq} a short exact sequence of the form
$$
\tau:\;0\to U_0\to U_i\to V^i\to 0.
$$
It follows from the long exact sequence
$$
0\to\Hom_A(V,U_0)\to\Hom_A(V,U_i)\to\Hom_A(V,V^i)\to\Ext^1_A(V,U_0)
$$
induced by $\tau$ that $\delta'_\tau(V)\leq j$.
Applying $(i-j)$ times a dual version of Lemma~\ref{cancelinseq} 
for the sequence $\tau$ we get that 
$U_i\simeq V^{i-j}\oplus\widetilde{U}_i$ for some module 
$\widetilde{U}_i$ of dimension
$$
d=\dim_kU_0+j\cdot\dim_kV.
$$
Hence the sequence $\sigma_{i+1}$ has the form
$$
0\to V^{i-j}\oplus\widetilde{U}_i\to V^{i+1-j}\oplus\widetilde{U}_{i+1}
\to V\to 0.
$$ 
We know that $\delta_{\sigma_{i+1}}(V)=0$, by Corollary~\ref{delta0to1}.
Then applying $(i-j)$ times Lemma~\ref{cancelinseq}
we get an exact sequence of the form
$$
0\to\widetilde{U}_i\to\widetilde{U}_{i+1}\oplus V\to V\to 0.
$$
Thus $\CO_{\widetilde{U}_i}\subseteq\ov{\CO}_{\widetilde{U}_{i+1}}$,
by Theorem~\ref{specialseq}.
Obviously the chain of inclusions 
$$
\ov{\CO}_{\widetilde{U}_j}\subseteq
\ov{\CO}_{\widetilde{U}_{j+1}}\subseteq
\ov{\CO}_{\widetilde{U}_{j+2}}\subseteq\ldots
$$
between orbit closures in $\mod_A(d)$ must terminate, i.e.,
$\ov{\CO}_{\widetilde{U}_l}=\ov{\CO}_{\widetilde{U}_{l+1}}$ for some $l>j$.
Hence $\widetilde{U}_l$ is isomorphic to $\widetilde{U}_{l+1}$ and
the sequence $\sigma_{l+1}$ splits.
\end{proof}

By our assumptions the sequence $\sigma_1$ does not split.
Moreover, the same holds for $\sigma_2$, by our choice of $g_1$.
Let $n\geq 1$ be the smallest integer such that the sequence
$\sigma_{n+2}$ splits.
As a direct consequence of Corollary~\ref{delta0to1} we get the following
fact.

\begin{cor} \label{deltai01}
$\delta_{\sigma_i}(V)=0$ and $\delta'_{\sigma_i}(V)=1$ for
$i\in\{1,2,\ldots,n+1\}$.
\end{cor}

Let $i\in\{1,2,\ldots,n\}$.
Since the sequence $\sigma_{i+1}$ does not split, $g_i$ does not factor
through the monomorphism $f_i$.
Consequently, the subspace $\CH_i$ of $\Hom_A(U_i,U_{i+1})$
spanned by $f_i$ and $g_i$ has dimension $2$.
The bilinear map
$$
F_i:\;\Ext^1_A(V,U_{i-1})\times\Hom_A(U_{i-1},U_i)\to\Ext^1_A(V,U_i),
$$
defined by $F_i(e,f)=\Ext^1_A(V,f)(e)$, induces a linear map
$$
F'_i:\;\Ext^1_A(V,U_{i-1})\to\Hom_k(\CH_i,\Ext^1_A(V,U_i)),
$$
such that $\left(F'_i(e)\right)(f)=F_i(e,f)=\Ext^1_A(V,f)(e)$.

\begin{lem} \label{injective}
The linear map $F'_i$ is injective for $i\leq n$.
\end{lem}

\begin{proof}
Let $\sigma:\;0\to U_{i-1}\to W\to V\to 0$ be a short exact sequence
such that $F'_i([\sigma])=0$.
Then the pushout of $\sigma$ under $f_i$ splits, and consequently,
$[\sigma]=c\cdot[\sigma_i]$ for some scalar $c\in k$, 
by Lemma~\ref{fromStoP} and Corollary~\ref{deltai01}.
We also know that the pushout of $\sigma$ under $g_i$ splits, which means
that
$$
0=\Ext^1_A(V,g_i)([\sigma])=c\cdot\Ext^1_A(V,g_i)([\sigma_i])
=c\cdot[\sigma_{i+1}].
$$
By our construction, $[\sigma_{i+1}]\neq 0$.
Thus $c=0$ and $[\sigma]=0$.
\end{proof}

Let $W$ and $W'$ be two vector spaces. If at least one of them is finite 
dimensional then the monomorphism
$$
F:\;W^\ast\otimes_kW'\to\Hom_k(W,W'),\qquad
F(\xi\otimes w')(w)=\xi(w)\cdot w',
$$
is an isomorphism.
Here $W^\ast$ denotes the dual space of $W$.
Thus $F'_i$ induces an injective linear map
$$
F''_i:\;\Ext^1_A(V,U_{i-1})\to\CH_i^\ast\otimes\Ext^1_A(V,U_i)
$$
for $i\in\{1,\ldots,n\}$.
Combining these maps we get an injective linear map
$$
F'':\;\Ext^1_A(V,U_0)\to\CH_1^\ast\otimes\CH_2^\ast\otimes\cdots\otimes
\CH_n^\ast\otimes\Ext^1_A(V,U_n).
$$

\begin{lem} \label{leq1}
$\rk(F'_1(e))=1$ for any $e\in\ov{\CE_{U_1}(V,U_0)}\setminus\{0\}$.
\end{lem}

\begin{proof}
We know that $\rk(F'_1(e))\geq 1$, by Lemma~\ref{injective}.
Assume now that $e$ belongs to $\CE_{U_1}(V,U_0)$.
It suffices to show that $\rk(F'_1(e))\leq 1$, as the function 
$\rk(F'_1(-)):\Ext^1_A(V,U_0)\to\BZ$ is semicontinuous.  
We take a short exact sequence
$$
\sigma:\;0\to U_0\to U_1\to V\to 0
$$
such that $[\sigma]=e$.
Since $\delta_\sigma(U_1)=1$, the kernel of the last map 
$\partial=\Ext^1_A(V,-)([\sigma])$ in the long exact sequence
$$
0\to\Hom_A(V,U_1)\to\Hom_A(U_1,U_1)\to\Hom_A(U_0,U_1)\xrightarrow{\partial}
\Ext^1_A(V,U_1)
$$
induced by $\sigma$ has codimension $1$.
We know that $\dim\CH_1=2$ and therefore there is a nonzero homomorphism
$h\in\CH_1$ such that 
$$
0=\partial(h)=\Ext^1_A(V,h)([\sigma])=F'_1([\sigma])(h).
$$
Consequently, $F'_1(e)$ is not injective and $\rk(F'_1(e))\leq 1$.
\end{proof}

If $e\in\ov{\CE_{U_1}(V,U_0)}\setminus\{0\}$ then $\Ker(F'_1(e))$ has dimension $1$.
Moreover,
$$
\Ker(F'_1(\lambda\cdot e))=\Ker(\lambda\cdot F'_1(e))=\Ker(F'_1(e))
$$
for any nonzero scalar $\lambda\in k$.
Thus we obtain a regular morphism
$$
\varphi_1:\BP\ov{\CE_{U_1}(V,U_0)}\to\BP\CH_1,\qquad
\varphi_1(k\cdot e)=\Ker(F'_1(e)).
$$

Let $\CH^\diamond_1=\CH_1\cap\CS_V(U_0,U_1)$.
It is an open, dense and $k^\ast$-invariant subset of $\CH_1$, 
by Lemma~\ref{Sopen}.
We can consider a regular map
$$
\widehat{\varphi}_1:\;\BP\CH^\diamond_1\to\BP\CE_{U_1}(V,U_0),\qquad
\widehat{\varphi}_1(k\cdot h)=\xi_{U_0,U_1,V}(h).
$$ 

\begin{lem} \label{immer}
$\varphi_1\circ\widehat{\varphi}_1(k\cdot h)=k\cdot h$ for any 
$h\in\CH^\diamond_1$.
\end{lem}

\begin{proof}
Let $h\in\CH^\diamond_1$. 
Then there is a short exact sequence of the form
$$
\sigma:\;0\to U_0\xrightarrow{h}U_1\to V\to 0.
$$
Obviously, $\widehat{\varphi}_1(k\cdot h)=k\cdot[\sigma]$.
Since the pushout of $\sigma$ under $h$ splits, $h$ belongs to the kernel
of $F'_1([\sigma])$.
But the kernel has dimension $1$ and therefore
$$
k\cdot h=\Ker(F'_1([\sigma]))=\varphi_1(k\cdot[\sigma])
=\varphi_1(\widehat{\varphi}_1(k\cdot h)).
$$
\end{proof}

\begin{lem} \label{singlefibre}
$\varphi_1^{-1}(k\cdot h)=\{\widehat{\varphi}_1(k\cdot h)\}$
for any $h\in\CH^\diamond_1$.
\end{lem}

\begin{proof}
Since $h\in\CS_V(U_0,U_1)$, there is a short exact sequence 
of the form
$$
\tau:\;0\to U_0\xrightarrow{h}U_1\to V\to 0.
$$
Let $\sigma:\;0\to U_0\xrightarrow{f'}W\to V\to 0$ be a short exact 
sequence such that $[\sigma]\in\ov{\CE_{U_1}(V,U_0)}\setminus\{0\}$ and
$\varphi_1(k\cdot[\sigma])=k\cdot h$.
We have to show that $k\cdot[\sigma]=\xi_{U_0,U_1,V}(h)$.
Since $h$ belongs to the kernel of $F'_i([\sigma])$, $h$ factors
through $f'$.
Hence $k\cdot[\sigma]=k\cdot[\tau]=\xi_{U_0,U_1,V}(h)$, 
by Lemma~\ref{fromStoP} and the equalities
$$
\delta_\tau(V)=\delta_{\sigma_1}(V)=0\qquad\text{and}\qquad
\delta'_\tau(V)=\delta'_{\sigma_1}(V)=1.
$$
\end{proof}

\begin{lem} \label{domin}
The map $\widehat{\varphi}_1:\;\BP\CH^\diamond_1\to\BP\CE_{U_1}(V,U_0)$ 
is dominant and the map $\varphi_1:\BP\ov{\CE_{U_1}(V,U_0)}\to\BP\CH_1$ 
is an isomorphism.
\end{lem}

\begin{proof}
Lemma~\ref{singlefibre} implies that 
$\widehat{\varphi}_1(\BP\CH^\diamond_1)=\varphi_1^{-1}(\BP\CH^\diamond_1)$
is an open subset of $\BP\ov{\CE_{U_1}(V,U_0)}$.
Moreover, $\BP\ov{\CE_{U_1}(V,U_0)}$ is irreducible, 
by Corollary~\ref{EMirreducible}.
Hence the map $\widehat{\varphi}_1$ is dominant and $\varphi_1$ is
birational.
Finally, $\varphi_1$ is an isomorphism as it is a projective map and
$\BP\CH_1$ is a smooth curve.
\end{proof}

We define inductively subsets $\CE_i\subseteq\Ext^1_A(V,U_{i-1})$ for
$i\in\{1,\ldots,n+1\}$ as follows:
$$
\CE_1=\ov{\CE_{U_1}(V,U_0)},\qquad\CE_{i+1}=F_i(\CE_i\times\CH_i).
$$
We consider the cofinite subset $\Lambda$ of $k$ consisting of the scalars
$\lambda$ such that $f_1+\lambda\cdot g_1$ belongs to $\CS_V(U_0,U_1)$,
or equivalently to $\CH^\diamond_1$. 
If $\lambda\in\Lambda$ then we derive from the diagram \eqref{diag1}
the following commutative diagram with exact rows:
\begin{equation} \label{diag2}
\vcenter{\xymatrix{
0\ar[r]&U_0\ar[rr]^{f_1+\lambda\cdot g_1}\ar[d]_{g_1}
 &&U_1\ar[r]\ar[d]^{g_2}&V\ar[r]\ar@{=}[d]&0\\
0\ar[r]&U_1\ar[rr]^{f_2+\lambda\cdot g_2}\ar[d]_{g_2}
 &&U_2\ar[r]\ar[d]^{g_3}&V\ar[r]\ar@{=}[d]&0\\
0\ar[r]&U_2\ar[rr]^{f_3+\lambda\cdot g_3}\ar[d]_{g_3}
 &&U_3\ar[r]\ar[d]^{g_4}&V\ar[r]\ar@{=}[d]&0.\\
&\vdots&&\vdots&\vdots
}}
\end{equation}
This implies the following fact.

\begin{cor} \label{aremono}
The homomorphism $f_i+\lambda\cdot g_i$ belongs to $\CS_V(U_{i-1},U_i)$ 
for any $\lambda\in\Lambda$ and $i\geq 1$. 
\end{cor}

We derive from Corollaries~\ref{deltai01} and \ref{aremono} 
a regular map 
$$
\eta_i:\Lambda\to\BP\CE_i,\qquad
\eta_i(\lambda)=\xi_{U_{i-1},U_i,V}(f_i+\lambda\cdot g_i),
$$
for any $i\in\{1,2,\ldots,n+1\}$.

\begin{lem} \label{twoinduction}
Let $i\in\{1,2,\ldots,n+1\}$. Then:
\begin{enumerate}
\item[(1)] The map $\eta_i:\Lambda\to\BP\CE_i$ is dominant.
\item[(2)] $\rk(F'_i(e))=1$ provided $i\leq n$ and $e\in\CE_i\setminus\{0\}$.
\end{enumerate}
\end{lem}

\begin{proof}
We proceed by induction on $i$.
For $i=1$ the claim follows from Lemmas~\ref{leq1} and \ref{domin}.
Assume now that the claim holds for some $i\leq n$.
By $(2)$, we can consider two regular maps
\begin{align*}
\varphi_i:\;&\BP\CE_i\to\BP\CH_i,&\varphi_i(k\cdot e)&=\Ker(F'_i(e)),\\
\psi_i:\;&\BP\CE_i\to\BP\CE_{i+1},&\psi_i(k\cdot e)&=\im(F'_i(e)).
\end{align*}
Moreover, $\psi_i$ is a surjective map such that
$$
\psi_i(\eta_i(\lambda))=\eta_{i+1}(\lambda),\quad\lambda\in\Lambda,
$$
by the diagram \eqref{diag2}.
Thus $\eta_{i+1}$ is a dominant map.
Assume now that $i<n$.
Observe that
$$
\CE^\diamond_{i+1}=\{e\in\CE_{i+1}\setminus\{0\};\;
k\cdot e\in\im(\eta_{i+1})\}
$$
is a dense subset of $\CE_{i+1}$. 
Hence by Lemma~\ref{injective}, it suffices to show that 
$\rk(F'_{i+1}(e))\leq 1$ for any $e\in\CE^\diamond_{i+1}$.
Let $e\in\CE^\diamond_{i+1}$.
Then there is a short exact sequence of the form
$$
\tau:\;0\to U_i\xrightarrow{f_{i+1}+\lambda\cdot g_{i+1}}U_{i+1}\to 
V\to 0
$$
with $[\tau]=e$, for some $\lambda\in\Lambda$.
Consequently $F'_{i+1}(e)(f_{i+1}+\lambda\cdot g_{i+1})=0$ and 
$F'_{i+1}(e)$ is not an injective map.
Thus $\rk(F'_{i+1}(e))<\dim_k\CH_{i+1}=2$.
\end{proof}

We get from the above proof the following regular maps:
$$
\xymatrix{
\BP\CE_1\ar[r]^{\psi_1}\ar[d]_{\varphi_1}
 &\BP\CE_2\ar[r]^{\psi_2}\ar[d]_{\varphi_2}
 &\cdots\ar[r]^{\psi_{n-1}}
 &\BP\CE_n\ar[r]^{\psi_n}\ar[d]_{\varphi_n}
 &\BP\CE_{n+1}\\
\BP\CH_1&\BP\CH_2&&\BP\CH_n.
}
$$

\begin{lem}
The map $\mu_i=\varphi_i\psi_{i-1}\ldots\psi_1:\BP\CE_1\to\BP\CH_i$
is an isomorphism for $i\in\{1,2,\ldots,n\}$.
\end{lem}

\begin{proof}
We consider the map $\eta=\mu_i\circ\varphi_1^{-1}:\BP\CH_1\to\BP\CH_i$
between projective lines.
It follows from the diagram \eqref{diag2} that 
$\eta(k\cdot(f_1+\lambda g_1))=k\cdot(f_i+\lambda g_i)$
for any $\lambda\in\Lambda$. 
Hence $\eta$ and $\eta\varphi_1=\mu_i$ are isomorphisms of varieties.
\end{proof}

\begin{lem} \label{En+1}
$\BP\CE_{n+1}=k\cdot[\sigma_{n+1}]$.
\end{lem}

\begin{proof}
Since the sequence $\sigma_{n+2}$ splits, the homomorphism $g_{n+1}$ 
factors through $f_{n+1}$.
By Lemma~\ref{twoinduction}(1) applied for $i=n+1$, it suffices to show
that $\xi_{U_n,U_{n+1},V}(f_{n+1}+\lambda g_{n+1})=k\cdot[\sigma_{n+1}]$ 
for any $\lambda\in\Lambda$.
The latter is true by Lemma~\ref{fromStoP}, Corollary~\ref{deltai01}
and since $f_{n+1}+\lambda g_{n+1}$ factors through $f_{n+1}$.
\end{proof}

Now we can prove the main result of the section.

\begin{prop} \label{keyprop}
$\BP\ov{\CE_{U_1}(V,U_0)}$ is a rational normal curve of degree $n$.
\end{prop}

\begin{proof}
We consider the image of $\BP\CE_1=\BP\ov{\CE_{U_1}(V,U_0)}$ 
via the linear map
$$
\BP F'':\BP\Ext^1_A(V,U_0)\to\BP\left(\CH_1^\ast\otimes\CH_2^\ast\otimes
\ldots\otimes\CH_n^\ast\otimes\Ext^1_A(V,U_n)\right).
$$
We know that $\CE_{n+1}$ is a linear subspace of dimension $1$,
by Lemma~\ref{En+1}.
Observe that
$$
\BP F''(\BP\CE_1)\subseteq\BP\left(\CH_1^\ast\otimes\CH_2^\ast\otimes
\ldots\otimes\CH_n^\ast\otimes\CE_{n+1}\right)
$$
and the latter space is canonically isomorphic to
$\BP\left(\CH_1^\ast\otimes\CH_2^\ast\otimes\ldots\otimes\CH_n^\ast\right)$.
Moreover, the induced map
$$
\CF:\BP\CE_1\to\BP\left(\CH_1^\ast\otimes\CH_2^\ast\otimes\ldots\otimes
\CH_n^\ast\right)
$$
is the composition of
$$
\BP\CE_1\xrightarrow{(\mu_1,\mu_2,\ldots,\mu_n)^T}
\BP\CH_1\times\BP\CH_2\times\ldots\times\BP\CH_n
$$
followed by the canonical isomorphism
$$
\BP\CH_1\times\BP\CH_2\times\ldots\times\BP\CH_n\to
\BP\CH_1^\ast\times\BP\CH_2^\ast\times\ldots\times\BP\CH_n^\ast
$$
and the Segre embedding 
$$
\BP\CH_1^\ast\times\BP\CH_2^\ast\times\ldots\times\BP\CH_n^\ast\to
\BP\left(\CH_1^\ast\otimes\CH_2^\ast\otimes\ldots\otimes\CH_n^\ast\right).
$$
Since the maps $\mu_i$ are isomorphisms, then in appropriate bases 
$\CF$ has the form
$$
(x:y)\mapsto\left(\ldots:x^{\varepsilon_1+\ldots+\varepsilon_n}
y^{n-\varepsilon_1-\ldots-\varepsilon_n}:\ldots\right)
_{(\varepsilon_1,\ldots,\varepsilon_n)\in\{0,1\}^n}.
$$
Thus $\CF(\BP\CE_1)$ is the rational normal curve of degree $n$ 
(see for instance \cite[Ex.1.14]{Harr}) and the same holds for 
$\BP\CE_1$.
\end{proof}

\section{Proof of Theorem~\ref{mainthm}}
\label{proof1}

Applying \cite[Thm.5]{Bmin} we get that there is a short exact sequence
$$
\sigma:\;0\to U\to M\to V\to 0
$$
and $\Sing(M,N)=\Sing(\ov{\CE_M(V,U)},0)$.
The assumption $\dim\CO_M-\dim\CO_N=2$ means that 
$$
2=[N,N]-[M,M]=\delta_\sigma(M)+\delta'_\sigma(N)
=\delta_\sigma(M)+\delta'_\sigma(V).
$$
If $\delta_\sigma(M)=0$ then $\Sing(M,N)=\Reg=\Sing(\CC_1,0)$, 
by \cite[Cor.3.7]{Zcodim2}.
Thus we may assume that $\delta_\sigma(M)\geq 1$.
We know that $\delta'_\sigma(V)\geq 1$, as the sequence $\sigma$ 
does not split.
Altogether, we get that $\delta_\sigma(M)=\delta'_\sigma(V)=1$.
Hence the claim follows from Proposition~\ref{keyprop}.
\qed

\section{Proof of Theorem~\ref{mainthm2}}
\label{proof2}

Let $(\Gamma_A,\tau_A)$ denote the Auslander-Reiten quiver of the algebra 
$A$ (see \cite{Rin} for details).
We may identify the set of vertices of $\Gamma_A$ with the corresponding
indecomposable modules.
A connected component $\CC$ of $\Gamma_A$ is called preprojective
if it has no cyclic paths and any $\tau_A$-orbit in $\CC$ contains
a projective module (see \cite{Rin}).
Moreover, a module is said to be preprojective if it is isomorphic to 
a direct sum of modules from preprojective components.

Let $M$ and $N$ be preprojective modules such that $\dim\CO_M-\dim\CO_N=2$
and $\CO_N\subseteq\ov{\CO}_M$.
By \cite[Thm.1.1]{Zcodim2}, we may assume that the modules $M$ and $N$ are 
disjoint.
Moreover we may assume that $N$ is a direct sum of at most two indecomposable
modules, by \cite[Thm.1.2]{Zcodim2}.
Repeating the proof of \cite[Lem.6.3]{Zcodim2}, we get that there are 
indecomposable direct summands $U$ and $V$ of $N$ such that
$$
[M,U]<[N,U],\;[U,M]=[U,N]
\quad\text{and}\quad
[M,V]=[N,V],\;[V,M]<[V,N].
$$
Therefore $N$ is isomorphic to $U\oplus V$ and the claim follows from
Theorem~\ref{mainthm}.
\qed

\section{Example}
\label{example}

Let $n\geq 3$ and $A$ be the path algebra of the quiver
$$
Q:\qquad
\vcenter
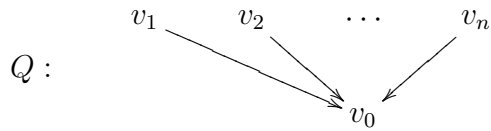
{\xymatrix{
v_1\ar[rrd]&v_2\ar[rd]&\cdots&v_n\ar[ld]\\
&&v_0
}}
$$
with $(n+1)$ vertices and $n$ arrows.
We may replace modules by representations of the quiver $Q$,
applying a Morita equivalence and \cite{Bgeo}.

Let $(a_1:b_1), (a_2:b_2),\ldots,(a_n:b_n)$ be pairwise different 
points of $\BP^1$. 
We define four representations of $Q$:
\begin{align*}
U:&\quad\vcenter{\xymatrix{
0\ar[rrd]&0\ar[rd]&\cdots&0\ar[ld]\\
&&k,
}}&
V:&\quad\vcenter{\xymatrix{
k\ar[rrd]_{[1]}&k\ar[rd]^(.67){[1]}&\cdots&k\ar[ld]^{[1]}\\
&&k,
}}\\
M:&\quad\vcenter{\xymatrix{
k\ar[rrd]_{\bsmatrix{a_1\\ b_1}}&k\ar[rd]^(.67){\bsmatrix{a_2\\ b_2}}
 &\cdots&k\ar[ld]^{\bsmatrix{a_n\\ b_n}}\\
&&k^2
}}
&\text{and}&\qquad N=U\oplus V.
\end{align*}
Then $\CS_V(U,M)\simeq\{(a,b)\in k^2;\;a_ib\neq b_ia\;\text{for}\;
i=1,2,\ldots,n\}$. 
In particular, $\CO_N\subseteq\ov{\CO}_M$, by Corollary~\ref{shortdeg}.
Since the representation $U$ is projective and the representation $V$ 
is injective, then
$$
[U,M]=[M,V]=2\qquad\text{and}\qquad[U,U]=[U,V]=[V,V]=1.
$$
Moreover, direct matrix calculations show that $[V,U]=0$ and $[M,M]=1$.
Hence
$$
[U,M]=[U,N],\qquad[M,V]=[N,V],\qquad[N,N]-[M,M]=2,
$$
and we can apply Theorem~\ref{mainthm}.
In fact, one can show that $\Sing(M,N)$ is the type of singularity of
a cone over a rational normal curve of degree $(n-2)$.
Observe that the representation $V$ is not preprojective 
(and $U$ is not preinjective) provided $n\geq 4$.


\bigskip

\noindent
Grzegorz Zwara\\
Faculty of Mathematics and Computer Science\\
Nicolaus Copernicus University\\
Chopina 12/18, 87-100 Toru\'n, Poland\\
E-mail: gzwara@mat.uni.torun.pl
\end{document}